\documentclass[twoside]{article}
\title{$(n,m)$-Strongly
Gorenstein Projective Modules}
\date{}
\usepackage[all]{xy}
\usepackage{amsfonts}
\usepackage{amsmath}
\usepackage{amssymb}
\usepackage{latexsym}
\newtheorem{thm}{\bf Theorem}[section]
\newtheorem{cor}[thm]{\bf Corollary}
\newtheorem{lem}[thm]{\bf Lemma}
\newtheorem{prop}[thm]{\bf Proposition}
\newtheorem{defn}[thm]{\bf Definition}

\catcode`\ç=13
\defç{\c{c}}
\catcode`\é=13
\defé{\'e}
\catcode`\à=13
\defà{\`a}
\catcode`\è=13
\defè{\`e}
\catcode`\â=13
\defâ{\^a}
\catcode`\ù=13
\defù{\`u}
\catcode`\ê=13
\defê{\^e}
\catcode`\î=13
\defî{\^\i}
\catcode`\ô=13
\defô{\^o}
\def\proof{{\parindent0pt {\bf Proof.\ }}}

\def\pd{{\rm pd}}
\def\fd{{\rm fd}}
\def\Gpd{{\rm Gpd}}

\def\lFFD{{ l\rm.FFD}}
\def\Im{{\rm Im}}
\def\Ker{{\rm Ker}}
\def\Ext{{\rm Ext}}

\def\Hom{{\rm Hom}}
\newcommand{\cqfd}
{\hspace{1cm}
\rule{2mm}{2mm}%
\medbreak%
\par%
}

\begin{document}
\thispagestyle{empty}
%%%%%%%%%%%%%%%%%%%%%%%%%%%%%%%%%%%%%%%%%%%%%%%%%%%%%%%%%
%%%%%%%%%%%%%%%%%%%%%%%%%%%%%%%%%%%%%%%%%%%%%%%%%%%%%%%%%
%%%%%%%%%%%%%%%%%%%%%%%%%%%%%%%%%%%%%%%%%%%%%%%%%%%%%%%%%
%%%TITLE%%%%%%%%%%%%%%%%%%%%%%%%%%%%%%%%%%%%%%%%%%%%%%%%%
\maketitle \vspace*{-2cm}
\begin{center}{\large\bf Driss Bennis}
%%%%%%%%%%%%%%%%%%%%%%%%%%%%%%%%%%%%%%%%%%%%%%%%%%%%%%%%%
%%%%%%%%%%%%%%%%%%%%%%%%%%%%%%%%%%%%%%%%%%%%%%%%%%%%%%%%%
%%%%%%%%%%%%%%%%%%%%%%%%%%%%%%%%%%%%%%%%%%%%%%%%%%%%%%%%%
%%%NAMES%%%%%%%%%%%%%%%%%%%%%%%%%%%%%%%%%%%%%%%%%%%%%%%%%

\bigskip
%%%%%%%%%%%%%%%%%%%%%%%%%%%%%%%%%%%%%%%%%%%%%%%%%%%%%%%%%
%%%%%%%%%%%%%%%%%%%%%%%%%%%%%%%%%%%%%%%%%%%%%%%%%%%%%%%%%
%%%%%%%%%%%%%%%%%%%%%%%%%%%%%%%%%%%%%%%%%%%%%%%%%%%%%%%%%
%%%%%%%%%%%%ADDRESSES%%%%%%%%%%%%%%%%%%%%%%%%%%%%%%%%%%%%%%%%%%%%%
\small{Department of Mathematics, Faculty of Science and
Technology of Fez,\\ Box 2202, University S. M.
Ben Abdellah Fez, Morocco, \\[0.12cm] driss\_bennis@hotmail.com}
\end{center}

\bigskip\bigskip
%%%%%%%%%%%%%%%%%%%%%%%%%%%%%%%%%%%%%%%%%%%%%%%%%%%%%%%%%
%%%%%%%%%%%%%%%%%%%%%%%%%%%%%%%%%%%%%%%%%%%%%%%%%%%%%%%%%
%%%%%%%%%%%%%%%%%%%%%%%%%%%%%%%%%%%%%%%%%%%%%%%%%%%%%%%%%
%%%ABSTRACT%%%%%%%%%%%%%%%%%%%%%%%%%%%%%%%%%%%%%%%%%%%%%%
\noindent{\large\bf Abstract.} This paper is a continuation of the
papers J. Pure Appl. Algebra, 210 (2007), 437--445 and J. Algebra
Appl., 8 (2009),  219--227. Namely, we introduce and study a
doubly filtered set of classes of  modules of finite Gorenstein
projective dimension, which are called $(n,m)$-strongly Gorenstein
projective ($(n,m)$-SG-projective for short) for integers $n\geq
1$ and $m\geq 0$. We are mainly interested in studying syzygies of
these modules. As consequences, we show that a module $M$ has
Gorenstein projective dimension at most $m$ if and only if
$M\oplus G$ is $(1,m)$-SG-projective for some Gorenstein
projective module $G$. And, over rings of finite left finitistic
flat dimension, that a module of finite Gorenstein projective
dimension has finite projective dimension if and only if it has
finite flat dimension.\bigskip

%%%%%%%%%%%%%%%%%%%%%%%%%%%%%%%%%%%%%%%%%%%%%%%%%%%%%%%%%
\medskip

%%%%%%%%%%%%%%%%%%%%%%%%%%%%%%%%%%%%%%%%%%%%%%%%%%%%%%%%%
%%%%%%%%%%%%%%%%%%%%%%%%%%%%%%%%%%%%%%%%%%%%%%%%%%%%%%%%%
%%%INTRODUCTION%%%%%%%%%%%%%%%%%%%%%%%%%%%%%%%%%%%%%%%%%%
\begin{section}{Introduction}
 Throughout this paper, $R$ denotes
a non-trivial associative ring with identity, and all modules are
left $R$-modules. For a module $M$, we use $\pd(M)$ and $\fd(M)$
to denote, respectively, the classical
projective and flat dimensions of $M$.\\
\indent  A module  $M$ is called \textit{Gorenstein projective}
(G-projective for short), if there exists an exact sequence of
projective modules,
$$\mathbf{P}=\ \cdots\rightarrow P_1\rightarrow P_0 \rightarrow
P_{ -1 }\rightarrow P_{-2 } \rightarrow\cdots,$$ such that  $M
\cong \Im(P_0 \rightarrow P_{ -1 })$ and such that $\Hom ( -, Q) $
leaves the sequence $\mathbf{P}$ exact whenever $Q$ is a
projective module. The exact sequence $\mathbf{P}$ is called a
complete projective resolution of  $M$.\\
\indent For a positive integer $n$, we say that $M$ has
\textit{Gorenstein projective dimension} at most $n$, and we write
$\Gpd_R(M)\leq n$ (or simply $\Gpd(M)\leq n$), if there is an
exact sequence of modules,
$$ 0 \rightarrow G_n\rightarrow \cdots \rightarrow G_0\rightarrow
M \rightarrow 0,$$ where each $G_i$ is Gorenstein projective
(suitable background materials on the notion of Gorenstein
projective modules can be found in \cite{LW,Rel-hom,HH}).\bigskip

The notion of Gorenstein projective modules was first introduced
and studied by  Enochs et al. \cite{GoInPj, GoIn, Fox} as a
generalization of the classical notion of projective  modules in
the sense that a module is projective if  and only if  it is
Gorenstein projective with finite projective dimension (see also
\cite{Rel-hom,HH}). In an unpublished work \cite[Theorem 4.2.6 and
Notes page 99]{LW},  Avramov, Buchweitz, Martsinkovsky, and Reiten
proved, over Noetherian rings, that finitely generated Gorenstein
projective modules  are just modules of Auslander's Gorenstein
dimension $0$ (\cite{A1}, see also \cite{A2}), which are
extensively studied by many others (part of the works on
Gorenstein dimension is summarized in Christensen's book
\cite{LW}).\bigskip

\indent The Gorenstein projective dimension has been extensively
studied by many others, who proved that this dimension shares many
nice properties of the classical projective dimension. In
\cite{BM}, Bennis and Mahdou introduced a particular case of
Gorenstein projective modules, which are defined as follows:

\begin{defn}[\cite{BM}]\label{defSG}
 \textnormal{A module $M$ is said to be strongly
Gorenstein projective (SG-projective for short), if there exists
an exact sequence of  projective modules, $$ \mathbf{P}=\
\cdots\stackrel{f}{\longrightarrow}P\stackrel{f}{\longrightarrow}P\stackrel{f}{\longrightarrow}P
\stackrel{f}{\longrightarrow}\cdots, $$ such that  $M \cong
\Im(f)$ and such that $\Hom ( -, Q) $ leaves the sequence
$\mathbf{P}$ exact whenever $Q$ is a projective module.}
\end{defn}
It is proved that the class of all strongly Gorenstein projective
modules is an intermediate class between the ones of projective
modules and Gorenstein projective modules \cite[Proposition
2.3]{BM}; i.e., we have the following inclusions
\begin{eqnarray*}
\{projective\ modules\}  &\subseteq&  \{SG\!-\!projective\ modules\}\\
  &\subseteq&  \{G\!-\!projective\ modules\}
\end{eqnarray*}
which are, in general, strict by \cite[Examples 2.5 and 2.13]{BM}.
The principal role of the strongly Gorenstein projective  modules
is to give the following characterization of Gorenstein projective
modules \cite[Theorem 2.7]{BM}: a module is Gorenstein projective
if and only if it is a direct summand of a strongly Gorenstein
projective module. The notion of strongly Gorenstein modules
confirm that there is an analogy between the notion of Gorenstein
projective modules and the notion of the usual projective modules.
In fact, this is obtained because the strongly Gorenstein
projective modules have simpler characterizations than their
correspondent Gorenstein  modules \cite[Propositions 2.9]{BM}. For
instance, a module $M$ is strongly Gorenstein projective if and
only if there exists a short exact sequence of modules,
$$0\rightarrow M\rightarrow P\rightarrow M\rightarrow 0,$$ where
$P$ is projective, and $\Ext(M,Q)=0$ for any projective module
$Q$. Using the results above, the notion of  strongly Gorenstein
projective modules was proven to be a good tool for establishing
results on Gorenstein projective dimension  (see, for instance,
\cite{BM2,BM3,BM4}). In \cite{BM2}, an extension of the notion of
strongly Gorenstein projective modules is introduced as follows:
for an integer $n>0$, a module $M$ is called $n$-strongly
Gorenstein projective ($n$-SG-projective for short), if there
exists an exact sequence of modules,
$$0\rightarrow M\rightarrow P_n\rightarrow\cdots\rightarrow P_1
\rightarrow M\rightarrow 0,$$ where each $P_i$ is projective, such
that $\Hom ( -, Q) $ leaves the sequence  exact whenever $Q$ is a
projective module (equivalently, $\Ext^{i}(M,Q)=0$ for $j+1\leq
i\leq j+n$ for some positive integer $j$ and for any projective
module $Q$ \cite[Theorem 2.8]{BM2}).  Then, 1-strongly Gorenstein
projective  modules are just strongly Gorenstein projective
modules. In \cite[Proposition 2.2]{BM2}, it is proved that an
n-strongly Gorenstein projective module is projective if and only
if it has finite flat dimension. In \cite{ZH}, Zhao and Huang,
continued the study of  $n$-strongly Gorenstein projective
modules. They gave more examples and they investigated the
relations between $n$-strongly Gorenstein projective modules and
$m$-strongly Gorenstein projective modules whenever $n\not= m$.
They  also proved, for two modules $M$ and $N$ projectively
equivalent (that is, there exist two projective modules $P$ and
$Q$ such  that $ M\oplus P\cong N\oplus Q$), that  $M$ is
$n$-strongly Gorenstein projective if and only if $N$ is
$n$-strongly Gorenstein projective  \cite[Theorem 3.14]{ZH} (see
Lemma \ref{lem-proj-equi} for a generalization of this result). So
using this result, we prove
the following lemma, which we use in the proof of the main results of this paper.\\
Recall, for a projective resolution of a module $M$, $$ \cdots
\rightarrow P_1 \rightarrow P_0\rightarrow M\rightarrow 0,$$ that
the module $K_i=\Im( P_i \rightarrow P_{i-1})$  for $i\geq 1$, is
called an $i^{th}$ syzygy of  $M$.

\begin{lem}\label{lem-stab-nSG} If $M$ is an $n$-strongly Gorenstein projective
module for some integer $n>0$, then:
\begin{enumerate}
\item Every $i^{th}$ syzygy of   $M$ is $n$-strongly Gorenstein projective.
    \item For every complete projective
resolution of $M$, $$\mathbf{P}=\ \cdots\rightarrow P_1\rightarrow
P_0 \rightarrow P_{ -1 }\rightarrow P_{-2 } \rightarrow\cdots,$$
every $\Im(P_i \rightarrow P_{i -1 })$ is $n$-strongly Gorenstein
projective.
\end{enumerate}
\end{lem}
\proof  First note that   $M$ admits a complete projective
resolution $$\mathbf{Q}=\ \cdots\rightarrow Q_2 \rightarrow Q_1
\rightarrow Q_{0}\rightarrow Q_{-1} \rightarrow\cdots   $$ in
which all images $\Im(Q_i\rightarrow Q_{i-1})$  are  $n$-strongly
Gorenstein projective modules. Indeed, $M$ is $n$-strongly
Gorenstein projective module, then there exists an exact sequence,
$$(*)\qquad 0\rightarrow M  \rightarrow
Q_{n -1}\rightarrow\cdots \rightarrow Q_0
 \rightarrow M\rightarrow 0,$$ where each
$Q_i$ is a projective module, such that $\Hom ( -, Q) $ leaves the
sequence exact whenever $Q$ is a projective module. For every
$i=1,...,n-1$, we decompose the exact sequence  $(*)$    into two
short exact sequences as follows: $$ \begin{array}{c}
 0\rightarrow M  \rightarrow Q_{n -1} \rightarrow\cdots \rightarrow Q_i \rightarrow N_i \rightarrow
0 \\\mathrm{and}\\  0\rightarrow N_i \rightarrow Q_{i-1}
\rightarrow\cdots  \rightarrow Q_0\rightarrow M\rightarrow 0
\end{array}
  $$
Assembling these sequences so that we obtain the following exact
sequence
$$0\rightarrow N_i\rightarrow Q_{i-1} \rightarrow \cdots\rightarrow
Q_0 \rightarrow Q_{n -1} \rightarrow \cdots \rightarrow  Q_i
\rightarrow N_i \rightarrow 0 $$ This shows that each $\Im(
Q_i\rightarrow Q_{i-1})$ is $n$-strongly Gorenstein projective.
Then, the desired complete  projective resolution $\mathbf{Q}$ is
obtained by assembling the sequence $(*)$ with
itself as done in the proof of \cite[Proposition 2.5(2)]{BM2}.\\
Now, using the left  half of $ \mathbf{Q}$, $ \cdots \rightarrow
Q_1 \rightarrow Q_0\rightarrow M\rightarrow 0 $, and the fact that
every two $i^{th}$ syzygies of a module are projectively
equivalent \cite[Theorem 9.4]{Rot}, the assertion $1$ follows from
\cite[Theorem 3.14]{ZH}.\\
We prove the second assertion. From $1$ it remains to prove the
result for the images of the right half of $\mathbf{P}$. Using
\cite[Proposition 1.8]{HH}, a dual proof of the one of
\cite[Theorem 9.4]{Rot} shows that the two module $\Im(Q_i
\rightarrow Q_{i-1})$ and $\Im(P_i \rightarrow P_{i-1})$ are
projectively equivalent for every $i \leq -1$, and therefore the
result follows from \cite[Theorem 3.14]{ZH}.\cqfd\bigskip

The aim of this paper is to generalize the notions above to a more
general context (Definition \ref{DefnSGproj}). Namely,  we
introduce and study a doubly filtered set of classes of  modules
with finite Gorenstein projective dimension, which are called
$(n,m)$-strongly  Gorenstein projective ($(n,m)$-SG-projective for
short) (for integers $n\geq 1$ and $m\geq 0$). First, we study the
relations between them (Proposition \ref{prop-relations}), and the
stability of this new class of modules under direct sum
(Proposition \ref{pro-sum}). Then, we set our first main result in
this paper (Theorem \ref{thm-main}), which shows, for an
$(n,m)$-SG-projective module $M$, that $\Gpd(M)=k\leq m$ for some
positive integer $k$. In particular, any $i^{th}$ syzygy  of $M$
is $(n,m-i)$-SG-projective for $1\leq i\leq k$, and any $i^{th}$
syzygy  of $M$ is $(n,0)$-SG-projective for $ i\geq k$. The second
main purpose of the paper is to investigate the converse of the
first main result. Namely, we ask: if an $i^{th}$ syzygy of a
module $M$ is $(n,m)$-SG-projective, is $M$  an
$(n,m+i)$-SG-projective module? In the second main result (Theorem
\ref{thm-converse-main}), we give an affirmative answer when $n=1$
as follows: for two integers $ d\geq 1$ and $m\geq 0$, if a
$d^{th}$ syzygy of a module $M$ is $(1,m)$-SG-projective, then
$\Gpd(M)=k\leq d+m$  for some positive integer $k$ and $M$ is
$(1,k)$-SG-projective. These results   lead to two results on
modules of finite
Gorenstein projective dimension:\\
\indent The first one shows that $(1,m)$-SG-projective modules can
serve to characterize modules of finite Gorenstein projective
dimension similarly to the characterization of Gorenstein
projective modules by strongly Gorenstein projective modules.
Namely, we prove (Corollary \ref{cor-char-gener}): for a module
$M$   and a positive integer $m$, $\Gpd(M)\leq m$ if and only if
$M\oplus G$ is $(1,m)$-SG-projective
for some Gorenstein projective module $G$.\\
\indent  The second one shows, over rings of finite left
finitistic flat dimension, that a module of finite Gorenstein
projective dimension has finite projective dimension if and only
if it has finite flat dimension (Proposition
\ref{prop-fd-pd-Gpd}). This, in fact, holds since we establish the
following extension of \cite[Proposition 2.2]{BM2} (Corollary
\ref{cor-fd-pd-n,m}): let $M$ be  an $(n,m)$-SG-projective module
for some integers $n\geq 1$ and $m\geq 0$. Then, $\pd(M)<\infty$
if and only if $\fd(M)<\infty$.

\end{section}
%%%%%%%%%%%%%%%%%%%%%%%%%%%%%%%%%%%%%%%%%%%%%%%%%%%%%%%%%
%%%%%%%%%%%%%%%%%%%%%%%%%%%%%%%%%%%%%%%%%%%%%%%%%%%%%%%%%
%%%%%%%%%%%%%%%%%%%%%%%%%%%%%%%%%%%%%%%%%%%%%%%%%%%%%%%%%
%%%%%%%%%%%%%%%%%%%%%%%%%%
%%%%%%%%%%%%%%%%%%%%%%%%%%%
%%%%%%%%%%%%%%%%%%%%%%%%%%%          SECTION 2:
%%%%%%%%%%%%%%%%%%%%%%%%%%%                       n-SG-Projective
%%%%%%%%%%%%%%%%%%%%%%%%%%%%%
%%%%%%%%%%%%%%%%%%%%%%%%%%%%%%%%%%%%%%%%%%%%%%%%%%%%%%%%%
%%%%%%%%%%%%%%%%%%%%%%%%%%%%%%%%%%%%%%%%%%%%%%%%%%%%%%%%%
%%%%%%%%%%%%%%%%%%%%%%%%%%%%%%%%%%%%%%%%%%%%%%%%
\begin{section}{Main results}  In this paper, we investigate the following kind of
modules:

\begin{defn}\label{DefnSGproj}\textnormal{Let $n\geq 1$ and $m\geq 0$ be
integers. A module $M$ is called $(n,m)$-SG-projective  if there
exists an exact sequence of modules, $$0\rightarrow M\rightarrow
Q_n\rightarrow\cdots\rightarrow Q_1 \rightarrow M\rightarrow 0,$$
where $\pd(Q_i)\leq m$  for $1\leq i\leq n $, such that $
\Ext^i(M,Q ) = 0 $ for any $i> m$   and for any projective module
$Q$.}
\end{defn}

Consequently, $(1,0)$-SG-projective modules are just strongly
Gorenstein projective  modules (by \cite[Proposition 2.9]{BM}),
and, generally, $(n,0)$-SG-projective modules are just
$n$-strongly Gorenstein projective  modules (by
\cite[Theorem 2.8]{BM2}).\\
One can show easily that modules of projective dimension at most
an integer $m$ are particular examples of  $(n,m)$-SG-projective
modules for every integer $n\geq 1$. The converse is not true in
general unless the $(n,m)$-SG-projective modules have finite flat
dimension  (see Corollary \ref{cor-fd-pd-n,m}). To give examples
of $(n,m)$-SG-projective modules with infinite projective
dimension, we can take any $(n,0)$-SG-projective module $M$ which
is not projective (use, for instance, \cite[Examples 2.4 and
2.6]{BM2} and \cite[Example 3.2]{ZH}) and any module $Q$ with
projective dimension at most $m$, then  we can show easily that
the direct sum $M \oplus Q$ is an $(n,m)$-SG-projective module
with infinite projective dimension.\bigskip

The main purpose of the paper is to investigate the syzygies of
$(n,m)$-SG-projective modules. In particular, we show that
$(n,m)$-SG-projective modules are particular examples of modules
with Gorenstein projective dimension at most $m$. Before, we give
some elementary properties of $(n,m)$-SG-projective modules.

\begin{prop}\label{prop-relations} Let $M$ be a module  and consider
two integers  $n\geq 1$ and $m\geq 0$.  We have the following
assertions:\begin{enumerate}
    \item If $M$ is $(n,m)$-SG-projective, then it is
    $(n,m')$-SG-projective for every $m'\geq m$.
    \item  If $M$ is $(n,m)$-SG-projective, then it is
    $(nk,m )$-SG-projective for every $k\geq 1$.\\
    In particular, every  $(1,m)$-SG-projective module  is
    $(n,m)$-SG-projective for every  $n\geq 1$.
\end{enumerate}
\end{prop}
\proof  1. Obvious.\\
2. Since    $M$ is $(n,m)$-SG-projective, there exists an exact
sequence of modules $0\rightarrow M\rightarrow
Q_n\rightarrow\cdots\rightarrow Q_1 \rightarrow M\rightarrow 0,$
where $\pd(Q_i)\leq m$  for $1\leq i\leq n $, such that $
\Ext^i(M,Q ) = 0 $ for any $i> m$ and for any projective module
$Q$. Assembling this sequence with itself $k$ times, we can show
that $M$ is also $(nk,m )$-SG-projective.\cqfd

\begin{prop}\label{pro-sum} Let  $(M_i)_{i\in I}$ be  a family of
 modules and consider the bounded families of integers $(n_i\geq 1)_{i\in I}$ and
$(m_i\geq 0)_{i\in I}$.\\ If,  for any $i\in I$, $M_i$ is
$(n_i,m_i)$-SG-projective, then the direct sum $\oplus_i \; M_i$
is $(n ,m )$-SG-projective, where $m=\max\{m_i\}$ and $n$ is the
least common multiple of $n_i$ for $ i\in I$.
\end{prop}
\proof  First, note that $m $  and $n$ exist since the families
$(n_i)_i$ and $(m_i)_i$ are bounded. Now, from Proposition
\ref{prop-relations}, $M_i$ is $(n ,m )$-SG-projective for any
$i\in I$. Then, using  standard arguments, we can show that the
direct sum $\oplus_i \; M_i$ is $(n ,m
)$-SG-projective.\cqfd\bigskip

Note, by \cite[Example 3.13]{ZH},  that the family of $(n ,m
)$-SG-projective modules is not closed under direct summands.
However, in Lemma \ref{lem-proj-equi} given later, we give a
situation in which a direct summand of an $(n,m)$-SG-projective
module is $(n,m)$-SG-projective.\bigskip

Now we give our first main result, in which we study the  syzygies
of an $(n,m)$-SG-projective.\\
Recall, for a projective resolution of a module $M$, $$\cdots
\longrightarrow P_1 \longrightarrow P_0\longrightarrow
M\longrightarrow 0,$$ that the module $K_i=\Im(P_i\rightarrow
P_{i-1})$ for $i\geq 1$ is called an $i^{th}$ syzygy of  $M$.

\begin{thm}\label{thm-main} If a module $M$ is
$(n,m)$-SG-projective for some integers $n\geq 1$ and $m\geq 0$,
then:\begin{enumerate}
    \item $\Gpd(M)=k\leq m$ for some positive integer $k$;
    \item Any $i^{th}$ syzygy  $K_i$ of $M$ is
    $(n,m-i)$-SG-projective for $1\leq i\leq k$;
    \item Any $i^{th}$ syzygy $K_i$  of $M$ is
    $(n,0)$-SG-projective for $ i\geq k$.
\end{enumerate}
\end{thm}
\proof  $1$ and $2$. Since $M$ is $(n,m)$-SG-projective, there
exists an exact sequence of modules,
$$(*)\quad 0\rightarrow M\rightarrow Q_n\rightarrow\cdots\rightarrow Q_1
\rightarrow M\rightarrow 0,$$ where $\pd(Q_i)\leq m$  for $1\leq
i\leq n $, such that $ \Ext^i(M,Q ) = 0 $ for any $i> m$ and for
any projective module $Q$. Consider a short exact sequence of
modules $$0 \rightarrow K_1 \rightarrow P_0\rightarrow
M\rightarrow 0,$$ where $P_0$ is projective. We  prove that $K_1$
is $(n,m-1)$-SG-projective. First, from \cite[Theorem 9.4]{Rot}, $
\Ext^i(K_1,Q ) = 0 $ for any $i> m-1$ and for any projective
module $Q$. Then, it remains to prove the existence of the exact
sequence. For that, decompose the exact sequence above $(*)$ into
short exact sequences
$$0\rightarrow H_{i}\rightarrow      Q_i\rightarrow
H_{i-1}\rightarrow 0,$$ where $H_{n}= M=H_0$ and
$H_{i}=\Ker(Q_{i}\rightarrow H_{i-1})$ for $i=1,...,n-1$. And
consider, for $i=0,...,n$, a short exact sequence $$0\rightarrow
K_{i,1} \rightarrow      P_{i,0}\rightarrow  H_{i}\rightarrow 0,$$
where $P_{i,0}$ is projective for $i=1,...,n-1$,  and
$P_{n,0}=P_{0,0}=P_0$,  and $K_{n,1}=K_{0,1}=K_1$. Applying  the
Horseshoe Lemma \cite[Lemma 6.20]{Rot}, we get the following
diagram for $i=n,...,1$:
$$\begin{array}{cccccccccc}
   &  &0&  &0&   &0&  &  \\
              &  &\uparrow  &  &\uparrow&   &\uparrow&  & \\
   0&\rightarrow &H_{i }&\rightarrow &Q_i& \rightarrow &H_{i-1}&\rightarrow &0 \\
              &  &\uparrow&  &\uparrow&   &\uparrow&  & \\
 0&\rightarrow &P_{i,0}&\rightarrow &P_{i,0}\oplus P_{i-1,0} & \rightarrow &P_{i-1,0}&\rightarrow &0 \\
                &  &\uparrow&  &\uparrow&   &\uparrow&  & \\
          0&\rightarrow &K_{i,1}&\rightarrow &Q'_i& \rightarrow &K_{i-1,1}&\rightarrow &0 \\
           &  &\uparrow&  &\uparrow&   &\uparrow&  & \\
   &  &0&  &0&   &0&  & \\
\end{array}$$
Assembling these diagrams we get the following diagram :
$$\begin{array}{ccccccccccccc}
    &  &0&  &0& & &  &0&   &0&  & \\
    &  &\uparrow&  &\uparrow& & &  &\uparrow&   &\uparrow&  & \\
   0&\rightarrow &M&\rightarrow &Q_n&\rightarrow&\cdots&\rightarrow &Q_1& \rightarrow &M&\rightarrow &0 \\
    &  &\uparrow&  &\uparrow& & &  &\uparrow&   &\uparrow&  & \\
   0&\rightarrow & P_0&\rightarrow &P_0\oplus P_{n-1,0}&\rightarrow&\cdots&\rightarrow &P_{1,0}\oplus P_0& \rightarrow &P_0&\rightarrow &0 \\
       &  &\uparrow&  &\uparrow& & &  &\uparrow&   &\uparrow&  & \\
      0&\rightarrow &K_1&\rightarrow &Q'_n&\rightarrow&\cdots&\rightarrow &Q'_1& \rightarrow &K_1&\rightarrow &0 \\
          &  &\uparrow&  &\uparrow& & &  &\uparrow&   &\uparrow&  & \\
     &  &0&  &0& & &  &0&   &0&  &
\end{array}$$
It is easy to show that $\pd(Q'_i)\leq m-1$  for $1\leq i\leq n $.
Hence, the bottom exact sequence of the diagram is the desired
sequence. Therefore, $K_1$ is $(n,m-1)$-SG-projective.\\
Then, by induction and using the same arguments above, we get that
$K_i$ is $(n,m-i)$-SG-projective for $i=1,...,m$. Particularly,
$K_m$ is $(n,0)$-SG-projective, then Gorenstein projective (from
\cite[Proposition 2.5]{BM2}), and
so $\Gpd(M)=k\leq m$ for some positive integer $k$.\\[0.2cm]
$ 3 $. Now, we prove that any $i^{th}$  syzygy of $M$ is
$(n,0)$-SG-projective for $ i\geq k$. Consider first $K_k$: a
$k^{th}$ syzygy of $M$. Since $K_k$ is Gorenstein projective, we
can chose a projective resolution of $K_k$ as a left half of any
of its complete projective resolution,  and so we get an exact
sequence
$$0 \rightarrow K'_{m-k}\rightarrow F_{m-k-1}\rightarrow \cdots\rightarrow
F_1\rightarrow F_0 \rightarrow  K_k \rightarrow 0,$$  where
$K'_{m-k}=\Im(F_{m-k}  \rightarrow F_{m-k-1}) $, such that $\Hom (
-, Q) $ leaves this sequence exact whenever $Q$ is a projective
module. From the first part of the proof,  $K'_{m-k}$ is
$(n,0)$-SG-projective (since it is an $m^{th}$ syzygy of $M$).
Then, dually to the first part of the proof, the dual version of
the Horseshoe Lemma \cite[Lemma 1.7]{HH} gives a raise to an exact
sequence of modules of the form: $$0\rightarrow  K_k \rightarrow
L_n\rightarrow\cdots\rightarrow L_1 \rightarrow K_k\rightarrow
0,$$ where $L_i$ is projective  for $1\leq i\leq n $. Then, with
the fact that $ \Ext^i(K_k,Q ) = 0 $ for any $i> 0$ and for any
projective module $Q$ (since $K_k$ is Gorenstein projective and by
\cite[Proposition 2.3]{HH}), we deduce that $K_k$ is
$(n,0)$-SG-projective. Therefore, from Lemma \ref{lem-stab-nSG}
with \cite[Theorem 3.14]{ZH}, we show that any $i^{th}$ syzygy
$K_i$ of $M$ is  $(n,0)$-SG-projective for $ i\geq
k$.\cqfd\bigskip

It is natural to ask for the converse of Theorem \ref{thm-main}.
Namely, we ask: if an $i^{th}$  syzygy of a module $M$ is
$(n,m)$-SG-projective, is  $M$  an   $(n,m+i)$-SG-projective
module? In the second main result, we give an affirmative answer
when $n=1$. For that, we need  the following two lemmas, which
are of independent interest.\\
\indent The first one gives a situation in which a direct summand
of an $(n,m)$-SG-projective module is $(n,m)$-SG-projective.

\begin{lem}\label{lem-proj-equi} Let $M$ and $N$ be two modules
such that $M\oplus P \cong N\oplus Q$ for some modules $P$ and $Q$
with finite projective dimension. Then, for two integers  $n\geq
1$ and $m\geq  \max\{\pd(P), \pd(Q)\}$, $M$ is   $(n ,m
)$-SG-projective if and only if $N$ is  $(n ,m )$-SG-projective.
\end{lem}
\proof  By symmetry, we only need to prove the direct implication.
The proof is   analogous to the one of \cite[Theorem 3.14]{ZH}.
For completeness, we give a proof here.\\ Since $M$ is   $(n ,m
)$-SG-projective, the direct sum $M\oplus P\cong N\oplus Q  $ is
also  $(n ,m )$-SG-projective (by Proposition \ref{pro-sum}).
Then, there exists for  $H=N\oplus Q $ an exact sequence of
modules,
$$0\rightarrow H\rightarrow Q_n\rightarrow\cdots\rightarrow Q_1
\rightarrow H\rightarrow 0 ,$$  where
 $\pd(Q_i)\leq m$  for $1\leq i\leq n $,
such that $ \Ext^i(H,L ) = 0 $ for any $i> m$   and for any
projective module $L$. Then, from \cite[Theorem 7.13]{Rot}, $
\Ext^i(N,L ) = 0 $ for any $i> m$   and for any projective module
$L$. Now, we have to construct the   exact sequence associated to
$N$. Decomposing the above sequence into three exact sequences:
 $$
\begin{array}{c}
   0\rightarrow H\rightarrow
Q_n\rightarrow E \rightarrow 0,\\ \quad 0\rightarrow E \rightarrow
Q_{n-1}\rightarrow \cdots\rightarrow Q_2 \rightarrow F \rightarrow
0,\quad \mathrm{and}  \\
 0 \rightarrow F\rightarrow Q_1 \rightarrow
H\rightarrow 0  \\
\end{array}  $$ Using the first and the last short exact
sequences above with, respectively, the trivial sequences
$0\rightarrow Q\rightarrow H\rightarrow N\rightarrow 0$ and
$0\rightarrow N\rightarrow H\rightarrow Q\rightarrow 0$, we get,
respectively, the following   pushout and pullback diagrams:\\
  \small{$\xymatrix{
     &    &  0 \ar[d] & 0 \ar[d]  &  \\
0\ar[r]& Q \ar@{=}[d] \ar[r] &H\ar[d] \ar[r] &N \ar@{-->}[d]
\ar[r] &
0\\
0  \ar[r]& Q \ar[r] & Q_n\ar[d] \ar@{-->}[r] & G_n \ar[d] \ar[r] &
0   \\
 &   & E\ar[d] \ar@{=}[r] & E\ar[d]   &
 \\
 &   & 0 &0  & } $  $\xymatrix{ \\ \\
\mathrm{and}   \\
 \\
 }$  $ \xymatrix{
     &  0 \ar[d]  & 0 \ar[d]  &   &  \\
 &  F\ar[d] \ar@{=}[r]&  F\ar[d]  &   &  \\
 0\ar[r]& G_1\ar@{-->}[d] \ar@{-->}[r] & Q_1\ar[d] \ar[r] & Q \ar@{=}[d]  \ar[r] & 0\\
0\ar[r]&N \ar[r] \ar[d]& H\ar[d] \ar[r]& Q\ar[r] & 0\\
 & 0 &0  &   & }$}\\
From Theorem \ref{thm-main}, $H$, $E$, and $F$ have Gorenstein
projective dimensions at most $m$. Then, from the diagrams above,
$G_1$ and $G_n$ have  finite  Gorenstein projective dimensions
which are, by standard  arguments, at most $m$. But, from the
middle sequence of each diagram, $G_1$ and $G_n$ have finite
projective dimensions. Then, from \cite[Proposition 2.27]{HH},
$\pd(G_1)=\Gpd(G_1)\leq m$ and $\pd(G_n)=\Gpd(G_n)\leq m$.
Finally, assembling the exact sequences: $$\begin{array}{c}
  0\rightarrow
N\rightarrow G_n\rightarrow E \rightarrow 0,\\ \quad 0\rightarrow
E \rightarrow Q_{n-1}\rightarrow \cdots\rightarrow Q_2 \rightarrow
F \rightarrow 0
,\quad \mathrm{and}  \\
0 \rightarrow F\rightarrow G_1
\rightarrow N\rightarrow 0 \\
\end{array}  $$  we get the following exact sequence:
$$0\rightarrow N\rightarrow G_n\rightarrow
Q_{n-1}\rightarrow\cdots\rightarrow Q_2 \rightarrow G_1
\rightarrow N\rightarrow 0 .$$ This completes the proof.\cqfd

\begin{lem} \label{lem-main2} Let $M$ be a module and let $n\geq 1$ and $m\geq 0$ be integers. Then,
\begin{enumerate}
    \item If $M$  is both Gorenstein projective and $(n,m)$-SG-projective,
then it is $(n,0)$-SG-projective.
    \item  If a $d^{th}$ syzygy of $M$ is
$(n,m)$-SG-projective  (for $ d\geq 1$), then $\Gpd(M)=k\leq d+m$
for some positive integer $k$ and any $i^{th}$ syzygy $K_i$  of
$M$ is $(n,0)$-SG-projective for $ i\geq k$.
\end{enumerate}
\end{lem}
\proof  $1.$ The proof is analogous to the last part of the proof
of Theorem \ref{thm-main}.\\[0.2cm]
$2.$ Since a $d^{th}$ syzygy of $M$ is $(n,m)$-SG-projective, we
can show  that $\Gpd(M)=k\leq d+m$  for some positive integer $k$.
Then, there exists an exact sequence of modules, $$ 0\rightarrow
K_k  \rightarrow P_{k-1}\rightarrow\cdots \rightarrow
P_0\rightarrow M\rightarrow 0,$$ where $P_i$ is projective for
$i=0,..., k-1$, and $K_k$ is Gorenstein projective.  Consider a
projective resolution of $K_k$ which is extracted from a left half
of one of its complete projective resolutions: $$  0\rightarrow
K_d \rightarrow Q_{d-1}\rightarrow\cdots Q_{k+1}\rightarrow
Q_{k}\rightarrow K_k \rightarrow 0,$$ where $Q_{k+i}$ is
projective for $i=0,..., d-k-1$, and $K_d$ is Gorenstein
projective. Clearly, $K_d$ is a $d^{th}$ syzygy of $M$. Hence, by
hypothesis, Lemma \ref{lem-proj-equi},  and since any two $i^{th}$
syzygies  of  $M$ are projectively equivalent, $K_d$  is
$(n,m)$-SG-projective, and then, from $(1)$, it is
$(n,0)$-SG-projective. This implies, by Lemma \ref{lem-stab-nSG},
that every  $\Im(Q_i \rightarrow Q_{i-1})$ is
$(n,0)$-SG-projective for $i\geq k+1$. Therefore, from Lemma
\ref{lem-proj-equi}, any $i^{th}$ syzygy $K_i$  of $M$ is
$(n,0)$-SG-projective for $ i\geq k$.\cqfd\bigskip

Now, we can prove the second main result:

\begin{thm}\label{thm-converse-main} Consider two integers $ d\geq 1$ and $m\geq 0$. If a $d^{th}$ syzygy of a module $M$ is
$(1,m)$-SG-projective, then $\Gpd(M)=k\leq d+m$  for some positive
integer $k$ and $M$ is $(1,k)$-SG-projective.
\end{thm}
\proof  By Lemma \ref{lem-main2} $(2)$, $\Gpd(M)=k\leq d+m$  for
some positive integer $k$ and any $i^{th}$ syzygy $K_i$  of $M$ is
$(1,0)$-SG-projective for $ i\geq k$. In particular, we have an
exact sequence of modules,
 $$  0\rightarrow K_k  \rightarrow
P_{k-1}\rightarrow\cdots \rightarrow P_0\rightarrow M\rightarrow
0,$$ where  $P_i$ is projective for $i=0,..., k-1$, and  the
$k^{th}$ syzygy $K_k$ of $M$  is $(1,0)$-SG-projective. Then,
there exists an exact sequence of modules, $$0\rightarrow
K_k\rightarrow P\rightarrow K_k\rightarrow 0,$$ where $P$ is
projective.  Then, by  \cite[Proposition 2.5(1) and its
proof]{BM2}, $K_k$ is $(k,0)$-SG-projective such that, by
assembling the short exact sequence above with itself $k$ times,
we have an exact sequence of the form $0\rightarrow K_k\rightarrow
P\rightarrow \cdots P\rightarrow K_k\rightarrow 0$. Then, using
the same proof as the one of \cite[Theorem 2.10]{HH}, we get the
following exact sequence: $$  0\rightarrow Q_k  \rightarrow
Q_{k-1}\rightarrow\cdots Q_1\rightarrow G\rightarrow M\rightarrow
0,$$ where $ Q_k =P$,  $ Q_i =P\oplus P_{i-1} $ for $i=1,...,
k-1$, and $ G =K_k\oplus P_0$. The module $ G =K_k\oplus P_0$ is
$(1,0)$-SG-projective with a short exact sequence $0\rightarrow
G\rightarrow Q\rightarrow G\rightarrow 0$, where $Q=P\oplus
P_0\oplus P_0$. Then, from the Horseshoe Lemma \cite[Lemma
6.20]{Rot}, we get the following diagram:
$$\begin{array}{cccccccccc}
   &  &0&  & &   &0&  &  \\
                 &  &\uparrow &  & &   &\uparrow&  & \\
   0&\rightarrow &M&  & &   &M&\rightarrow &0 \\
              &  &\uparrow&  & &   &\uparrow&  & \\
   0&\rightarrow &G&\rightarrow &Q& \rightarrow &G&\rightarrow &0 \\
              &  &\uparrow&  &\uparrow&   &\uparrow&  & \\
 0&\rightarrow &Q_1&\rightarrow &Q_1\oplus Q_1 & \rightarrow &Q_1&\rightarrow &0 \\
                &  &\uparrow&  &\uparrow&   &\uparrow&  & \\
                &  &    \vdots&  &\vdots&   &\vdots&  & \\
            &  &\uparrow&  &\uparrow&   &\uparrow&  & \\
          0&\rightarrow &Q_{k-1}&\rightarrow &Q_{k-1} \oplus Q_{k-1}& \rightarrow &Q_{k-1}&\rightarrow &0 \\
           &  &\uparrow&  &\uparrow&   &\uparrow&  & \\
            0&\rightarrow &Q_{k}&\rightarrow &Q'_{k}& \rightarrow &Q_{k}&\rightarrow &0 \\
                       &  &\uparrow&  &\uparrow&   &\uparrow&  & \\
   &  &0&  &0&   &0&  & \\
\end{array}$$
Since $Q_{k}$ is projective, $Q'_{k}$ is projective. Putting the
cokernel  into this diagram, we obtain an exact sequence
$0\rightarrow M\rightarrow P' \rightarrow M\rightarrow 0$ such
that, by the middle exact sequence, $\pd(P')\leq k$. Therefore,
$M$ is $(1,k)$-SG-projective.\cqfd\bigskip

As consequences of the two main results, we get some results on
modules with finite Gorenstein projective dimension.\bigskip

The first one extends the role of strongly Gorenstein projective
modules (i.e., $(1,0)$-SG-projective modules), which serve to
characterize Gorenstein projective modules, to the setting of
$(1,m)$-SG-projective modules as follows:

\begin{cor}\label{cor-char-gener}
Let $M$ be a module and let $m$ be a positive integer. Then,
$\Gpd(M)\leq m$ if and only if there exists a Gorenstein
projective module $G$ such that the direct sum $M\oplus G$ is
$(1,m)$-SG-projective.
\end{cor}
\proof   $\Leftarrow.$ Follows from Theorem \ref{thm-main}(1) and
\cite[Proposition 2.19]{HH}.\\[0.2cm]
$\Rightarrow.$ Since $\Gpd(M)\leq m$, there exists an exact
sequence of modules, $$(*)\quad 0 \rightarrow K_m\rightarrow
P_{m-1}\rightarrow\cdots \rightarrow P_0\rightarrow M \rightarrow
0,$$ where $P_i$ is projective for $i=0,..., m-1$ and $K_m$ is
Gorenstein projective. Then, from \cite[Theorem 2.7]{BM}, there
exists a Gorenstein projective module $G'$ such that  $K_m\oplus
G'$ is $(1,0)$-SG-projective. From the right half of a complete
projective resolution of $G'$, we get an exact sequence, $$ 0
\rightarrow G'\rightarrow Q_{m-1}\rightarrow\cdots \rightarrow
Q_0\rightarrow G \rightarrow 0,$$  where $Q_i$ is projective for
$i=0,..., m-1$ and $G$ is Gorenstein projective. Adding this
sequence with the sequence $(*)$, we get the following exact
sequence  $$ 0 \rightarrow K_m\oplus G'\rightarrow P_{m-1} \oplus
Q_{m-1}\rightarrow\cdots \rightarrow P_{0} \oplus Q_0\rightarrow
M \oplus G \rightarrow 0.
$$  This means that the $m^{th}$ syzygy $K_m\oplus G'$ of $M
\oplus G$ is $(1,0)$-SG-projective. Therefore, from Theorem
\ref{thm-converse-main}, $M\oplus G$ is
$(1,m)$-SG-projective.\cqfd\bigskip

The second corollary investigates the relation between
$(n,m)$-SG-projective modules and the usual projective dimension.
It is known, for a module $M$, that $\Gpd(M)\leq \pd(M)$ with
equality if $\pd(M)<\infty$. For $(n,m)$-SG-projective modules we
have the following result, which is an extension of
\cite[Proposition 2.2]{BM2}:

\begin{cor}\label{cor-fd-pd-n,m}
Let $M$ be  an $(n,m)$-SG-projective module for some integers
$n\geq 1$ and $m\geq 0$. Then, $\pd(M)<\infty$ if and only if
$\fd(M)<\infty$.
\end{cor}
\proof    We only need to proof the converse implication. Assume
that $\fd(M)<\infty$, then so every syzygy of $M$ has finite flat
dimension. From Theorem \ref{thm-main}, an $m^{th}$ syzygy of $M$
is $(n,0)$-SG-projective, and so it is projective from
\cite[Proposition 2.2]{BM2}. This implies that $\pd(M)<\infty$, as
desired.\cqfd\bigskip

The above result leads us to conjecture that every module of
finite Gorenstein projective dimension has finite projective
dimension if it has finite flat dimension. From \cite[Corollary
2.3]{BM2}, we have an affirmative answer over rings with finite
weak global dimension. In the following result, we give an
affirmative answer in a more general context. Recall that the left
finitistic flat dimension of $R$ is the quantity
$\lFFD(R)=\sup\{\fd_R(M)\;|\;M\;is\;an \;
 R\!-\!module\;with$ $\;\fd_R(M)<\infty\} $.

\begin{prop}\label{prop-fd-pd-Gpd}
If $\lFFD(R)<\infty$, then  every module with both finite
Gorenstein projective dimension and finite flat dimension has
finite projective dimension.
\end{prop}
\proof  Assume that $\lFFD(R)=n$ for some positive integer $n$.
Let $M$ be a module such that $\fd(M)<\infty $ and $\Gpd(M)=k<
\infty$. To see that $\pd(M)<\infty $, it is sufficient, from
 Corollary \ref{cor-char-gener} and its proof, to show that
$K_m\oplus G'$ is projective (we use the notation of Corollary
\ref{cor-char-gener} and its proof). From the proof of
\cite[Theorem 2.7]{BM}, $K_m\oplus G'$ can be considered as the
direct sum of all the images of a complete projective resolution
of $K_m$. Now, since $\fd(K_m)\leq n$ (since $\fd(M)<\infty $),
all the images of  this complete projective resolution have finite
flat dimension, which is at most $n$ (since $\lFFD(R)=n$). This
implies that $\fd(K_m\oplus G')\leq n$. Therefore, from
\cite[Proposition 2.2]{BM2}, $\pd(K_m\oplus G')<\infty$, as
desired.\cqfd\bigskip

Finally, it is convenient to note that one could define and study
$(n,m)$-SG-injective modules  as a dual notion to the current one
of $(n,m)$-SG-projective modules. Then, every result established
here for  $(n,m)$-SG-projective modules, except Corollary
\ref{cor-fd-pd-n,m} and Proposition \ref{prop-fd-pd-Gpd}, has a
dual version for $(n,m)$-SG-injective modules.\bigskip

\noindent {\bf Acknowledgment.} The author would like to thank the
referee for his/her careful reading of this work.\bigskip

\end{section}

%%%%%%%%%%%%%%%%%%%%%%%%%%%%%%%%%%%%%%%%%%%%%%%%%%%%%%%%%
%%%REFERENCES%%%%%%%%%%%%%%%%%%%%%%%%%%%%%%%%%%%%%%%%%%%%
%%%%%%%%%%%%%%%%%%%%%%%%%%%%%%%%%%%%%%%%%%%%%%%%%%%%%%%%

%%%%%%%%%%%%%%%%%%%%%%%%%%%%%%%%%%%%%%%%%%%%%%%%%%%%%%%%
\end{document}